\documentclass[journal]{IEEEtran}  

\IEEEoverridecommandlockouts

\usepackage{graphicx}
\usepackage{epsfig} 
\usepackage{times} 
\usepackage{amsmath}
\usepackage{amssymb}  
\usepackage{color}
\usepackage{amsfonts}
\usepackage{subfigure}
\usepackage{multirow}
\usepackage{multicol}
\usepackage{siunitx}
\usepackage{overpic}
\usepackage{mathrsfs}  

\usepackage{booktabs}
\usepackage{threeparttable}

\usepackage{epstopdf}
\usepackage{xspace}

\begin{document}

\title{\LARGE \bf On 1-D PDE-Based Cardiovascular Flow Bottleneck Modeling and Analysis: A Vehicular Traffic Flow-Inspired Approach}

\author{Nikolaos Bekiaris-Liberis

\thanks{N. Bekiaris-Liberis is with the Department of Electrical \& Computer Engineering, Technical University of Crete, Chania, Greece, 73100. Email address: \texttt{bekiaris-liberis@ece.tuc.gr}.}}

\maketitle
\pagestyle{headings}

\markboth{Submitted to IEEE Transactions on Automatic Control on October 22, 2021}{Bekiaris-Liberis}

\begin{abstract}
We illustrate the potential of PDE-based \textcolor{black}{traffic flow} control in cardiovascular flow analysis, monitoring, and control, presenting a PDE-based control-oriented formulation, for 1-D blood flow dynamics in the presence of stenosis. This is achieved adopting an approach for modeling and analysis that relies on the potential correspondence of 1-D blood flow dynamics in the presence of stenosis, with 1-D traffic flow dynamics in the presence of bottleneck. We reveal such correspondence in relation to the respective (for the two flow types), speed dynamics and a (consistent with them) fundamental diagram-based reduction; bottleneck dynamic effects description and resulting boundary conditions; and free-flow/congested regimes characterization.

\end{abstract}
\section{Introduction}
\label{secintro}

Arterial stenosis, due to, for example, atherosclerotic plaque building up in arteries or in-stent restenosis, is a primary cause of human losses worldwide \cite{Takahata_1}. A great number of deceases, attributed to congested blood flow, currently accounting for about 50\% of deaths within the European Union \cite{Quarteroni_0}, could be avoided with accurate/timely detection and action implementation. This is true particularly in view of the practical feasibility that is supported by existing technologies, such as, for example, smart, stents and bypass grafts, and other implantable devices, where actuation and sensing may be performed wirelessly, via communication with a central computer; see, for example, \cite{Takahata_1}, \cite{takahata ena}, \cite{kiourti}, \cite{Takahata_0}, \cite{takahata dyo}. 



Despite the technological advancement and urgent need for availability of respective advanced methodologies, illustrated by
their potential in congested blood flow detection/treatment, \textcolor{black}{there exists no {\em control-theoretic} approach tackling the formulation, analysis, monitoring, and feedback control problems of congested blood flow}, in its natural, continuous in time/space, domain, in the presence of stenosis. However (and despite the domain and dimensional complexity of cardiovascular flow), there exist 1-D, second-order, hyperbolic Partial Differential Equation (PDE) systems that may effectively describe (on average) blood flow dynamics; see, e.g., \cite{Canic0}, \cite{Formaggia_2}, \cite{Formaggia_1}, \cite{Li-Canic}, \cite{Quarteroni_1}, \cite{Wang}. Thus, such models may be utilized as basis for control-theoretic modeling, analysis, and design purposes.


In this paper we launch an effort in this direction formulating and analyzing, from a \textcolor{black}{\em PDE-based} {\em traffic flow control} (see, e.g., \cite{Blandin}, \cite{Claudel}, \cite{Maria Laura}, \cite{Goatin}, \cite{karafyllis}, \cite{Piacentini}, \cite{Tumash 1}, \cite{Tumash 2}, \cite{Yu}, \cite{Prieur}) perspective, the dynamics of 1-D blood flow in the presence of stenosis. The stenosis is considered to be located at the boundary of the arterial segment considered. We present two alternative formulations in which the stenosis dynamics are characterized via a static or dynamic description for the pressure drop through the stenosis. Together with utilization of a baseline dynamic model for blood flow, capturing the main transport phenomena and respective mass/momentum conservation principles, such formulation gives rise to a $2\times 2$ (heterodirectional; see, e.g., \cite{Auriol}, \cite{Rafa}) hyperbolic PDE system, with a static or dynamic boundary condition, at the outlet of the artery segment considered, respectively. As the location, geometry, and length of the potential stenosis are considered to be unknown, the derived model may incorporate unknown PDE domain length and boundary conditions parameters. 


For the derived dynamic descriptions of 1-D blood flow in the presence of stenosis, we then illustrate the correspondence, of certain features, with traffic flow dynamics in the presence of bottleneck. We explore correspondence with Payne-Whitham- and Aw-Rascle-Zhang-type models, in particular, in relation to speed dynamics and a consistent (with respect to reduction to conservation law equation, for instance, of Lighthill-Whitham-Richards-type) fundamental diagram. Furthermore, we illustrate the connection to respective, dynamic models of traffic flow bottlenecks. In particular, boundary blood flow stenosis may be characterized via the pressure drop at the stenosis location, while boundary traffic flow bottleneck may be described via the capacity drop at the bottleneck area. Moreover, for each type of stenosis description, we provide a consistent boundary condition at the outlet, which could either be static or dynamic, also illustrating the correspondence with the respective boundary conditions, in the case of traffic flow bottleneck. We also discuss the analogy between characterization of free-flow/congested traffic regimes and supercritical/subcritical blood flow regimes.





We start in Section \ref{sec model} presenting a control-oriented model for blood flow in which arterial stenosis is described either as static or dynamic, boundary bottleneck. In Section \ref{sec analysis 1} we analyze the obtained hyperbolic system, revealing the dynamic correspondence with traffic flow dynamics in the presence of bottlenecks. In Section \ref{conlus} we discuss potential research directions that may emerge from the results presented.

\section{Control-Theoretic Modeling of Stenosis}
\label{sec model}
\subsection{Baseline 1-D cardiovascular flow model}
We consider the following $2\times 2$ hyperbolic system, which constitutes an 1-D approximation of cardiovascular flow dynamis (see, e.g., \cite{Formaggia_2}, \cite{Formaggia_1}, \cite{Quarteroni_1})
\setlength{\arraycolsep}{0pt}\begin{eqnarray}
A_t(x,t)&=&-A_x(x,t)V(x,t)-A(x,t)V_x(x,t)\label{sys1}\\
V_t(x,t)&=&-V(x,t)V_x(x,t)-\frac{1}{\rho}\frac{\partial P\left(A(x,t)\right)}{\partial x}\nonumber\\
&&-K_r\frac{V(x,t)}{A(x,t)}\label{speed}\\
A(0,t)V(0,t)&=&Q_{\rm in}(t)\label{boundary qin},
\end{eqnarray}\setlength{\arraycolsep}{5pt}where
\noindent
$A>0$ is section area of artery, $V>0$ is average blood speed, $\rho>0$ is blood density, $K_r>0$ is friction parameter related to blood viscosity, $t\geq 0$ is time, $x\in[0,D]$ is spatial variable, $D>0$ is length of artery segment considered, $P\in\mathbb{R}$ is pressure, which accounts for vessel wall displacement, and $Q_{\rm in}>0$ is flow at the inlet of the artery segment considered (it could, for example, be described by a periodic signal, with period equal to the cardiac cycle, see, e.g.,~\cite{Quarteroni_1}). A possible expression for the pressure function is (see, e.g., \cite{Formaggia_1}, \cite{Quarteroni_1})
\begin{eqnarray}
P\left(A\right)&=&\frac{\beta}{A_0}\left(\sqrt{A}-\sqrt{A_0}\right)\label{pressure}\\
\beta&=&{hE}\sqrt{\pi}b,\label{beta_def}
\end{eqnarray}
where $A_0$ is reference arterial section area at rest, $h>0$ is artery wall thickness, $E>0$ is Young's modulus, and $b$ is a positive parameter. One boundary condition, associated with (\ref{sys1}), (\ref{speed}), is provided in (\ref{boundary qin}), describing the blood flow entering the arterial segment considered. The second boundary condition is specified in the following sections since it depends both on the sign of the eigenvalues of hyperbolic system (\ref{sys1}), (\ref{speed}) as well as the stenosis dynamic description adopted.

In the present setup, in which the case of a boundary bottleneck is considered, it is assumed that \textcolor{black}{parameters $\beta$ and $A_0$, in the pressure equation (\ref{pressure}), are known and constant} throughout the domain, which may be a reasonable requirement given that variations in geometry and mechanical properties of the artery, imposed by the stenosis, are considered to be located at the boundary $x=D$. Although most of modeling and analysis developments could be performed considering spatially-varying coefficients $\beta$, $A_0$, for presentation and formulation simplicity, as well as to not distracting the reader from the main scope of the paper, which is presentation and analysis of a control-theoretic, stenosis model and its correspondence with traffic flow bottleneck model, we do not consider this case here.

\subsection{Formulation of available measurements output equation}

In the present paper we consider the case in which the pressure and flow at the inlet of the artery segment considered are measured in real time. Although such a setup may appear, at first sight, as unrealistic, current technological advancements enable the availability of these measurements. In particular, such measurements could be wirelessly transmitted to a central computer, utilizing smart stent (or bypass graft) devices, see, for example, \cite{Takahata_1}, \cite{kiourti}, \cite{Takahata_0}, \cite{takahata dyo}. Thus, besides having available $Q_{\rm in}$, a measured output is available, given by
\begin{eqnarray}
y(t)=P\left(A(0,t)\right).
\end{eqnarray}
Since location, geometry, and material properties of the stent, in realistic scenarios, could be considered as known, it follows that $\beta$ and ${A_0}$ at $x=0$ are known (even in the case in which $\beta$, $A_0$ may take different values, as compared with their values for $x\in(0,D)$). Thus, using (\ref{pressure}), measurements of $A(0,t)$ could be obtained, and hence, using (\ref{boundary qin}), measurements of $V(0,t)$.

\subsection{Stenosis model as static boundary bottleneck}
\label{boundary stenosis}

\begin{figure}[h]
\centering
\includegraphics[width=\linewidth]{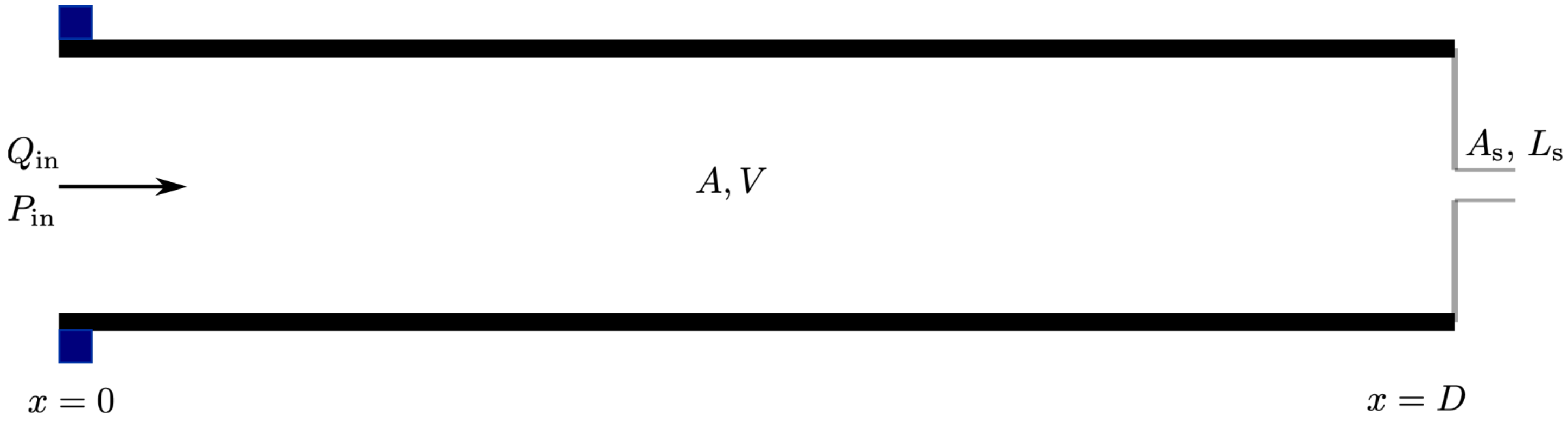}
\caption{Simplified schematic of an 1-D approximation of an arterial segment \textcolor{black}{at rest} with boundary stenosis. At the inlet of the segment considered (i.e., for $x=0$) an implantable smart stent device measures pressure and flow (i.e., $P_{\rm in}$ and $Q_{\rm in}$, respectively). The bottleneck location $D$, stenosis cross-sectional area $A_{\rm s}$, and stenosis length $L_{\rm s}$ may be unknown.}
\label{fig2}
\end{figure}

This potential formulation of a bottleneck is derived assuming that the stenosis (e.g., due to atherosclerotic plaque building up at arterial wall \cite{Quarteroni_0} or in-stent restenosis \cite{Takahata_1}) is located downstream of the inlet (i.e., the known location of a, for instance, stent device). In particular, we treat the right boundary of the arterial segment considered as the point at which the potential stenosis is located. Therefore, the spatial variable $x$ belongs to $\left[0,D\right]$, where $D$ may be unknown as the stenosis location may be unknown. The right boundary condition is derived such that it incorporates the effect of stenosis in the outlet. A schematic view of the setup considered is shown in Fig. \ref{fig2}.


\paragraph{Modeling assumptions} In the present setup, the domain length $D$ and the effective section area at the stenosis location $A_{\rm s}$ may be unknown\footnote{The stenosis may be assumed to be axisymmetric (see, for example, \cite{young-tsai}; see also, for instance, \cite{Lee}, \cite{Misra}, for more complex, potential stenosis geometries).}. In particular, the stenosis section area is assumed to be constant, which may be a reasonable assumption given the material and elastic properties of atherosclerotic plaque (see, for example, \cite{Takashima}). It is further assumed that (for constant $A_{\rm s}$) flow is conserved through the stenosis. \textcolor{black}{Moreover, the least complex formulation of the stenosis effect (at least in terms of the number of potentially unknown model parameters) could be obtained assuming (initially) zero (or, effectively, very small) length for the stenosis.}


\paragraph{Boundary condition formulation} Consequently, at the stenosis location the following relation may be satisfied (\cite{seely}, \cite{young-tsai0}, \cite{young-tsai}; see also \cite{clark}, \cite{Koeppl}, \cite{stergiopulos} for relevant expressions)
\begin{eqnarray}
\!\!\Delta P\left(A(D,t),V(D,t)\right)\!\!\!\!\!\!&=&\!\!\!\!\!\!V(D,t)^2\frac{K_{\rm s}\rho}{2}\left(\frac{A(D,t)}{A_{\rm s}}-1\right)^2\!\!,\label{stenosis no length}
\end{eqnarray}
where $\Delta P\left(A(D),V(D)\right)$ denotes the pressure drop due to the stenosis, while parameter $K_{\rm s}>0$ is known (obtained, for example, from experimental data, see, e.g., \cite{seely}). The pressure drop denotes the pressure difference between the locations before and after the stenosis. For the former we may assume that it is given by (\ref{pressure}), while for the latter, we may assume that it is described such that a terminal boundary condition, modeling the effect of blood flow dynamics in arteries downstream of the stenosis (see, for example, \cite{Formagga_3}, \cite{stergiopulos}), is imposed. Therefore, we may define 
\begin{eqnarray}
\!\!\!\!\!\!\!\Delta P\left(A(D,t),V(D,t)\right)\!\!\!\!\!\!&=&\!\!\!\!\!\!P\left(A(D,t)\right)\!-\!R_{\rm T}A(D,t)V(D,t),\label{drop1}
\end{eqnarray}
with $Q_{\rm s}=A(D)V(D)$ denoting the flow at the inlet of the stenosis, where
$R_{\rm T}\geq0$ denotes a total, terminal resistance. Parameter $R_{\rm T}$ may be chosen depending on the blood flow conditions modeled for a considered arterial network, and thus, it may be considered as known. Using (\ref{pressure}), (\ref{stenosis no length}), (\ref{drop1}) we obtain
\begin{eqnarray}
&&\frac{\beta}{A_0}\left(\sqrt{A(D,t)}-\sqrt{A_0}\right)-R_{\rm T}A(D,t)V(D,t)\nonumber\\
&&-V(D,t)^2\frac{K_{\rm s}\rho}{2}\left(\frac{A(D,t)}{A_{\rm s}}-1\right)^2=0\label{stenosis no length11}.
\end{eqnarray}
Equation (\ref{stenosis no length11}) prescribes a boundary condition at $x=D$, associated with system (\ref{sys1})--(\ref{boundary qin}), with $A_{\rm s}$ and $D$ being unknown.

Although, as starting point and under the assumption of zero stenosis length, formulation (\ref{stenosis no length11}) may appear to be adequately realistic, a more accurate, nevertheless more complex, formulation for the right boundary condition (at the inlet of the stenosis) may be obtained utilizing the following relation for the pressure drop (see, for example, \cite{seely}, \cite{young-tsai0}, \cite{young-tsai})
\setlength{\arraycolsep}{0pt}\begin{eqnarray}
\Delta P\left(A(D,t),V(D,t)\right)&=&\frac{8{\pi}\mu L_{\rm s}}{A_{\rm s}^2}A(D,t)V(D,t)+V(D,t)^2\nonumber\\
&&\times\frac{K_{\rm s}\rho}{2}\left(\frac{A(D,t)}{A_{\rm s}}-1\right)^2,\label{stenosis length}
\end{eqnarray}\setlength{\arraycolsep}{5pt}where $L_{\rm s}>0$ is unknown stenosis length and $\mu$ is known blood viscosity coefficient. Thus, using (\ref{pressure}), (\ref{drop1}), (\ref{stenosis length}) we obtain 
\begin{eqnarray}
&&\frac{\beta}{A_0}\left(\sqrt{A(D,t)}-\sqrt{A_0}\right)-\left(R_{\rm T}+\frac{8{\pi}\mu L_{\rm s}}{A_{\rm s}^2}\right)A(D,t)\nonumber\\
&&\times V(D,t)-V(D,t)^2\frac{K_{\rm s}\rho}{2}\left(\frac{A(D,t)}{A_{\rm s}}-1\right)^2=0\label{stenosis no length11new}.
\end{eqnarray}
In the case of boundary condition (\ref{stenosis no length11new}), in addition to $A_{\rm s}$ and $D$, the stenosis length $L_{\rm s}$ may also be an unknown parameter.

\subsection{Stenosis model as dynamic boundary bottleneck}
\label{boundary stenosis dynamic}
Boundary condition formulations (\ref{stenosis no length11}), (\ref{stenosis no length11new}) may be accurate for zero or, effectively, very small, stenosis length. A potentially more realistic, nevertheless more complex, model of the pressure drop dynamics, accounting for larger stenosis length (yet, much smaller than the length $D$ of the arterial segment considered), may be written as (see, for example, \cite{stergiopulos}, \cite{young-tsai})
\setlength{\arraycolsep}{0pt}\begin{eqnarray}
V_t(D,t)&=&\frac{1}{\rho L_{\rm s}}\Delta P\left(A(D,t),V(D,t)\right)-V(D,t)^2\frac{K_{\rm s}}{2L_{\rm s}}\nonumber\\
&&\!\!\times\left(\frac{A(D,t)}{A_{\rm s}}-1\right)^2-\frac{8{\pi}\mu }{\rho A_{\rm s}^2} A(D,t)V(D,t).\label{dod}
\end{eqnarray}\setlength{\arraycolsep}{5pt}Employing (\ref{pressure}), (\ref{drop1}), relation (\ref{dod}) may be written as
\setlength{\arraycolsep}{0pt}\begin{eqnarray}
\!\!{V}_t(D,t)&=&\frac{\beta}{\rho A_0 L_{\rm s}}\left(\sqrt{A(D,t)}-\sqrt{A_0}\right)-\left(\frac{A(D,t)}{A_{\rm s}}-1\right)^2\nonumber\\
&&\!\!\!\!\!\!\!\!\!\times \frac{K_{\rm s}}{2L_{\rm s}} V(D,t)^2-\left(\frac{8{\pi}\mu }{\rho A_{\rm s}^2}+\frac{R_{\rm T}}{\rho L_{\rm s}} \right)A(D,t)V(D,t).\label{ODE-boundary1}
\end{eqnarray}\setlength{\arraycolsep}{5pt}The complete model (\ref{sys1})--(\ref{boundary qin}), (\ref{ODE-boundary1}) consists of a nonlinear hyperbolic PDE - nonlinear Ordinary Differential Equation (ODE) coupled system.

We note here that, similarly to considering an ODE for describing the dynamics at the right boundary, due to the presence of stenosis, one could consider an ODE for the dynamics of the pressure downstream of the stenosis (instead of the static relation given by the term $R_{\rm T}AV$ in (\ref{drop1}); see, for example, \cite{stergiopulos}). For formulation and presentation simplicity, we do not consider the dynamics of the stenosis pressure here.



\section{Analysis of the cardiovascular flow model and its relation to traffic flow dynamics}
\label{sec analysis 1}

\subsection{Analysis of the hyperbolic system}
\label{section 2ac}

\paragraph{Blood flow information propagation}\textcolor{black}{In physiological} conditions blood flow is reported to lie in congested (or, subcritical) regime (see, e.g., \cite{Formaggia_1}, \cite{Quarteroni_1}). In particular, the eigenvalues of the hyperbolic system (\ref{sys1}), (\ref{speed}) are given by
\begin{eqnarray}
\bar{\lambda}_1\left(A,V\right)&=&V+\sqrt{\frac{\beta}{2\rho A_0}}A^{\frac{1}{4}}\label{eigen1}\\
\bar{\lambda}_2\left(A,V\right)&=&V-\sqrt{\frac{\beta}{2\rho A_0}}A^{\frac{1}{4}}.\label{eigen2}
\end{eqnarray}

Since we are concerned with the case of subcritical regime we restrict our attention in a nonempty, connected open subset $\Omega$ of the set $\bar{\Omega}=\left\{\left(A,V\right)\in\mathbb{R}^2:0<A,0<V\right\}$, such that $V<\sqrt{\frac{\beta}{2\rho A_0}}A^{\frac{1}{4}}$, and hence, $\bar{\lambda}_1>0$ and $\bar{\lambda}_2<0$, in the region of interest. System (\ref{sys1})--(\ref{boundary qin}) is then strictly hyperbolic with distinct, real nonzero eigenvalues, as long as $\left(A,V\right)\in{\Omega}$, which implies that information propagates both forward (with blood flow) and backward (at a lower speed).

\paragraph{Transformation to Riemann variables}
\label{riemann sub}
The Riemann invariants that correspond to the hyperbolic system (\ref{sys1}), (\ref{speed}) with eigenvalues (\ref{eigen1}), (\ref{eigen2}) are defined as
\begin{eqnarray}
u\left(A,V\right)&=&V+2\sqrt{\frac{2\beta}{\rho A_0}}A^{\frac{1}{4}}\label{Reiman1}\\
v\left(A,V\right)&=&V-2\sqrt{\frac{2\beta}{\rho A_0}}A^{\frac{1}{4}}.\label{Reiman2}
\end{eqnarray}
The inverse transformations that correspond to (\ref{Reiman1}), (\ref{Reiman2}) are 
\begin{eqnarray}
V\left(u,v\right)&=&\frac{1}{2}\left(u+v\right)\label{vreim}\\
A\left(u,v\right)&=&\frac{\rho^2A_0^2}{4^5\beta^2}\left(u-v\right)^4.
\end{eqnarray}
In the new variables, system (\ref{sys1}), (\ref{speed}) is written as
\begin{eqnarray}
u_t(x,t)&=&-{\lambda}_1\left(u(x,t),v(x,t)\right)u_x(x,t)\nonumber\\
&&+f_1\left(u(x,t),v(x,t)\right)\label{transformed system}\\
v_t(x,t)&=&-{\lambda}_2\left(u(x,t),v(x,t)\right)v_x(x,t)\nonumber\\
&&+f_1\left(u(x,t),v(x,t)\right),\label{transformed system1}
\end{eqnarray}
where
\setlength{\arraycolsep}{0pt}\begin{eqnarray}
{\lambda}_1\left(u,v\right)&=&\frac{5u+3v}{8}\\
{\lambda}_2\left(u,v\right)&=&\frac{3u+5v}{8}\\
f_1\left(u,v\right)&=&-\frac{4^{\frac{9}{2}}K_{r}\beta^2}{\rho^2A_0^2}\frac{u+v}{\left(u-v\right)^4}.\label{deff1}
\end{eqnarray}\setlength{\arraycolsep}{5pt}Boundary condition (\ref{boundary qin}) at the inlet is expressed in terms of the Reimann variables as 
\begin{eqnarray}
g\left(u(0,t),v(0,t)\right)&=&Q_{\rm in}(t)\label{boundary qin1}\\
g\left(u,v\right)&=&\frac{\rho^2A_0^2}{4^{\frac{11}{2}}\beta^2}\left(u+v\right)\left(u-v\right)^4.\label{g_func}
\end{eqnarray}
Together with (\ref{transformed system})--(\ref{g_func}) we associate a boundary condition at $x=D$, which may be specified as follows.

\paragraph{Boundary condition at the outlet}Since the $2\times 2$ hyperbolic system (\ref{transformed system}), (\ref{transformed system1}) is heterodirectional, together with the boundary condition (\ref{boundary qin1}) at $x=0$, one should specify a boundary condition at $x=D$. There are different options for specifying a boundary condition at $x=D$, also depending on the coupling type, of the arterial segment considered, with different arteries (also considering different types of arteries; see, e.g., \textcolor{black}{\cite{Formagga_3}, \cite{Koeppl}, \cite{Quarteroni_1}}). 

Since in the present paper we are concentrated on the modeling of bottleneck effects, the boundary condition is specified in order to describe the pressure difference between the locations before and after the bottleneck, also accounting for a cumulative effect of arteries downstream of the stenosis area. This could be achieved employing a static (see Section~\ref{boundary stenosis}) or dynamic (see Section \ref{boundary stenosis dynamic}) description for the effect of the stenosis. For completeness we also discuss the case in which there is no stenosis and arteries downstream of the arterial segment considered do not affect its dynamics. In such a case, we explore a non-reflecting (see, e.g., \cite{Formaggia_2}, \cite{thompson}) boundary condition at $x=D$, such that there is no incoming wave at the right boundary. We summarize below these three cases.

$\quad$1) In the case of static boundary bottleneck, the boundary condition at $x=D$ is specified in order to describe the pressure \textcolor{black}{drop} at the outlet of the arterial segment considered, where a stenosis is located, as described in Section~\ref{boundary stenosis}. Specifically, using (\ref{stenosis no length11new}), the right boundary condition is expressed in Riemann variables as
\begin{eqnarray}
G\left(u(D,t),v(D,t)\right)&=&0\label{boundary qin2},
\end{eqnarray}
where
\setlength{\arraycolsep}{0pt}\begin{eqnarray}
G\left(u,v\right)&=&\rho \left(u-v\right)^2-\frac{32\beta}{\sqrt{A_0}}-d_1\left(u-v\right)^4\left(u+v\right)\nonumber\\
&&-{4K_{\rm s}\rho}\left(u+v\right)^2\left(d_2\left(u-v\right)^4-1\right)^2\label{f_func}\\
d_1&=&\left(R_{\rm T}+\frac{8{\pi}\mu L_{\rm s}}{A_{\rm s}^2}\right)\frac{\rho^2A_0^2}{ \beta^2 4^3},\quad d_2\!=\!\frac{\rho^2A_0^2}{4^5\beta^2 A_{\rm s}}.
\end{eqnarray}\setlength{\arraycolsep}{5pt}

$\quad$2) Similarly, in the case of a dynamic boundary bottleneck, the right boundary condition is expressed in Riemann variables using (\ref{ODE-boundary1}), (\ref{vreim}) as
\begin{eqnarray}
v(D,t)&=&2X(t)-u(D,t)\label{boundary qin2dyn}\\
\dot{X}(t)&=&\frac{1}{32\rho L_{\rm s}}G\left(u(D,t),2X(t)-u(D,t)\right)\label{f_funcdyn}.
\end{eqnarray}

$\quad$3) In the case in which there is no stenosis, we explore the option of a non-reflecting, right boundary condition (see, e.g., \cite{Formaggia_2}, \cite{thompson}), such that no incoming wave enters at the right boundary of the arterial segment considered. Such a boundary condition could be compared with a ``free'' right boundary condition (imposed, for example, in specific traffic networks, see, e.g., \cite{bekiaris}, \cite{karafyllis}; see also Section \ref{subtra}). Such a boundary condition could be described as
\begin{eqnarray}
v_t(D,t)&=&f_1\left(u(D,t),v(D,t)\right).\label{dynamic bound}
\end{eqnarray}
In fact, one could observe that, boundary condition (\ref{dynamic bound}) implies (for classical solutions) that the Riemann variable corresponding to the negative eigenvalue has zero spatial derivative at the right boundary.

\textcolor{black}{Well-posedness} of the $2\times 2$ hyperbolic PDE system (\ref{transformed system})--(\ref{g_func}) with the dynamic boundary condition (\ref{dynamic bound}) (with (\ref{deff1})) or (\ref{boundary qin2dyn}), (\ref{f_funcdyn}) (with (\ref{f_func})), or the static boundary condition (\ref{boundary qin2}) (with (\ref{f_func})) may be guaranteed utilizing, for example, the results in \cite{bastin}, \cite{TI-Global}. To be able to employ such results, certain assumptions are required to be imposed on regularity, size, and compatibility (with boundary conditions) of initial conditions, on size and regularity of flow $Q_{\rm in}$ at the inlet, and on the values of parameters $\beta$, $A_0$, $K_{\rm s}$, $\rho$, $L_{\rm s}$, $\mu$, $R_{\rm T}$, $A_{\rm s}$. Well-posedness of the hyperbolic system considered, for realistic values of the various parameters involved, is also consistent with the dynamic behavior of the actual, physical system (see, for example, \cite{Canic0}, \cite{Quarteroni_1}). Although important, we do not belabor this issue as it is beyond the present paper's primary scope.  

\subsection{Properties of the model from a traffic flow perspective}
\label{subtra}
$\quad$1) The first correspondence with second-order traffic flow models originates in the speed equation (\ref{speed}). Such relation (for $K_{r}=0$) may be compared to the speed dynamics of Payne-Whitham traffic flow model (see, e.g., \cite{treiber}) with pressure given by (\ref{pressure}). Equation (\ref{sys1}), which expresses the conservation of blood volume entering and exiting an artery segment considered, corresponds to the conservation of the number of vehicles entering and exiting a given highway segment. 

$\quad$2) The correspondence of model (\ref{sys1}), (\ref{speed}) to traffic flow models of Payne-Whitham (and Aw-Rascle-Zhang, see, e.g., \cite{Fan}, \cite{lebaque}, \cite{treiber}) type could be also viewed via a fundamental diagram definition, considering the pressure function (\ref{pressure}), which could be explained as follows. Adopting the procedure in \cite{zhang} for derivation of a fundamental diagram relation from the speed equation (\ref{speed}), we define $V=F(A)$ and substitute this relation into (\ref{speed}) in order to obtain only one, conservation law equation of the form (\ref{sys1}), i.e., of the form $A_t+\left(F(A)A\right)_x=0$, where $F$ is to be determined. With $K_{r}=0$ we get that
\setlength{\arraycolsep}{0pt}\begin{eqnarray}
\!\!\!\!\!\!F'(A)\left(A_t+F(A)A_x+\frac{\beta}{2\rho A_0 \sqrt{A}F'(A)}A_x\right)&=&0.\label{see1}
\end{eqnarray} \setlength{\arraycolsep}{5pt}Imposing the reasonable requirement that $F'(A)<0$, for all $A>0$, relation (\ref{see1}) holds if the following equation is satisfied 
\begin{eqnarray}
A_t+F(A)A_x+\frac{\beta}{2\rho A_0 \sqrt{A}F'(A)}A_x&=&0.\label{see2}
\end{eqnarray} 
Therefore, in order for equation (\ref{see2}) to reduce to the conservation law equation (\ref{sys1}), imposing $V=F(A)$, for any $A$, the following should hold
\begin{eqnarray}
F'(A)^2&=&\frac{\beta}{2\rho A_0 }A^{-\frac{3}{2}}.
\end{eqnarray} 
Therefore, since $F'(A)<0$, for any $A>0$, we get that
\begin{eqnarray}
F(A)&=&F(0)-2\sqrt{\frac{{2}{\beta}}{{\rho A_0} }}A^{\frac{1}{4}}.\label{fd}
\end{eqnarray} 
The constant $F(0)$ may be viewed as the speed at a limiting case in which the section area tends to zero. Thus, in practice, it may be defined, for example, through considering a maximum possible, blood transport speed, which could be obtained empirically. Relation (\ref{fd}) defines a fundamental diagram (see, e.g., \cite{treiber}), satisfying the required conditions. In particular, function $\bar{Q}(A)=AF(A)=A\left(F(0)-2\sqrt{\frac{{2}{\beta}}{{\rho A_0} }}A^{\frac{1}{4}}\right)$, for $A\in[0,A_1]$, where  $A_1=\frac{F(0)^4\rho^2A_0^2}{\beta^2 4^3}$, satisfies $\bar{Q}(0)=\bar{Q}\left(A_1\right)=0$, while being strictly concave.  

We note here that the limiting case in which $V=F(A)$, constitutes a considerable simplification, which may appear, at first sight, as not realistic for cardiovascular systems. However, such a reduction may be useful in, for example, studying the dynamic effect of a bottleneck in blood flow, at a vicinity upstream of the stenosis, employing only the respective conservation law equation.

$\quad$3) The two different bottleneck descriptions also bear a resemblance to traffic flow bottleneck descriptions. For example, boundary bottlenecks may appear due to lane-drops or, in general, due to the presence of locations of reduced capacity, at the end of a controlled area of interest, such as, for example, where a tunnel or an area of high curvature begins (see, e.g., \cite{papageorgiou}). A boundary bottleneck could be described through properly modeling the traffic capacity drop at the bottleneck location (potentially also employing different fundamental diagram relations for the traffic speed immediately before and after the bottleneck location; see, e.g., \cite{treiber}, \cite{Yu-Koga}); corresponding to the static equation (\ref{stenosis no length}) (or (\ref{stenosis length})), which describes the pressure drop at the area of the stenosis (that may also be viewed as defining a pressure fundamental diagram at the stenosis, depending on the pressure immediately before, as $P_{\rm s}\left(A,V\right)=P\left(A\right)-V^2\frac{K_{\rm s} \rho}{2}\left(\frac{A}{A_{\rm s}}-1\right)^2$, which becomes a function of only $A$ when $V=F(A)$).

In the case of a dynamic bottleneck description, speed (or flow) dynamics at the area of the stenosis are described by an ODE (as in (\ref{ODE-boundary1}), (\ref{f_funcdyn}); see also, e.g., \cite{stergiopulos}, \cite{young-tsai}), dictated by the pressure difference between the areas at the inlet and outlet of the stenosis. This may be viewed as corresponding to the case of dynamic description of traffic density at a bottleneck area through an ODE, dictated by the flow difference between the flow arriving and exiting the bottleneck area (see, e.g., \cite{bekiaris-krstic}, \cite{papageorgiou}). In both cases the resulting dynamic description consists of a nonlinear, hyperbolic PDE-ODE coupled system.

$\quad$4) The non-reflecting boundary condition (see, e.g., \cite{Formaggia_2}, \cite{thompson}), considered in the case of no bottleneck, at the outlet of the arterial segment considered, could be compared with a free boundary condition in traffic flow, at the outlet of a highway segment (see, e.g., \cite{bekiaris}, \cite{karafyllis}). In the case of traffic flow, a free boundary condition indicates that, at the right boundary, there is no influence from the downstream traffic. Respectively, a non-reflecting boundary condition, in the case of blood flow, indicates that no incoming wave enters the domain at the right boundary. It could be verified that, in both cases, the respective dynamic boundary conditions imply a zero spatial derivative (for classical solutions) of the leftward transporting Riemann variable, at the right boundary.



$\quad$5) For cardiovascular flow, subcritical regime is characterized by the sign of $\bar{\lambda}_2$ in (\ref{eigen2}). Analogously, traffic congestion may be characterized by negative sign of a respective eigenvalue that corresponds to the Reimann invariant transporting opposite to traffic flow direction (see, e.g., \cite{belleti}). One difference lies in that physiological conditions for cardiovascular flow correspond to the subcritical (congested) regime (where $\bar{\lambda}_2<0$; see, e.g., \cite{Formaggia_1}, \cite{Quarteroni_1}), whereas for traffic flow, physiological conditions may be considered as corresponding to the free-flow (supercritical) regime (where $\bar{\lambda}_2>0$; see also, e.g., \cite{belleti}).

\section{Future Perspectives}
\label{conlus}
\textcolor{black}{The arterial stenosis models presented and the correspondence with vehicular traffic flow bottleneck models, may constitute the starting point for PDE-based, control-theoretic developments for cardiovascular flow stenosis analysis, estimation, and control, inspired by respective traffic flow techniques.}

A first research direction, of practical and theoretical significance, would be to develop algorithms for real-time identification of potential stenosis location and characteristics (such as, e.g., length and thickness). Towards this end, a possible approach would be to design adaptive observers, utilizing the $2\times 2$ hyperbolic system (\ref{transformed system}), (\ref{transformed system1}), (\ref{boundary qin1}) with either (\ref{boundary qin2}) or (\ref{boundary qin2dyn}), (\ref{f_funcdyn}), aiming at simultaneous state estimation and parameters identification, employing the derived model and the available boundary measurements. \textcolor{black}{Related methods for traffic flow models and general hyperbolic systems could be found, for example, in \cite{anfinsen}, \cite{bin}, \cite{seo}, \cite{ji wang}, \cite{Yu-estimation}.}

Another potential research direction would be to consider the feedback control problem of blood flow at areas with stenosis. Towards this end, perhaps a crucial issue would be to specify how, in practice, the required actuation could be performed. One possibility would be to consider boundary actuation, manipulating the inflow in (\ref{boundary qin1}) through certain micro-electromechanical systems (for example, smart stent devices, actuated wirelessly; see, e.g., \cite{takahata ena}, \cite{takahata dyo}), thus resulting in a boundary control problem for system (\ref{transformed system}), (\ref{transformed system1}), (\ref{boundary qin1}) with either (\ref{boundary qin2}) or (\ref{boundary qin2dyn}), (\ref{f_funcdyn}). An alternative possibility would be to consider in-domain actuation enabled through automated drug delivery systems (see, e.g., \cite{takahata ena}, \cite{kulkani}). Such an approach could build upon an extension of the presented model to incorporate in-domain actuation, in correspondence with automated vehicles-based actuation incorporation in vehicular traffic (see, for example, \cite{bekiaris}, \cite{darbha}, \cite{karafyllisd}, \cite{Piacentini}, \cite{yi-horowitz}).

\section*{Acknowledgments}
The author would like to thank Prof. Argiris Delis and Prof. Ioannis Nikolos for fruitful discussions.

\end{document}